\documentclass[conference]{IEEEtran}

\ifCLASSINFOpdf
\else
\fi

\usepackage[numbers,sort&compress]{natbib}
\usepackage{amsmath, amssymb}
\usepackage{graphicx}
\usepackage{float}
\usepackage[utf8]{inputenc}
\usepackage{mathtools}
\usepackage{cleveref}
\usepackage{url} 
\usepackage{enumerate} 
\usepackage[super]{nth}
\usepackage{multicol}

\usepackage{cellspace}

\usepackage{parskip}
\begin{document}
\setlength{\abovedisplayskip}{6pt}
\setlength{\belowdisplayskip}{6pt}	
	
\bstctlcite{IEEEexample:BSTcontrol}
 
\title{Self-Scheduling Robust Preview Controllers for Path Tracking and Autonomous Vehicles \\  }

\author
	{\IEEEauthorblockN{Ali Boyali\IEEEauthorrefmark{1},
		Lyu Zheming, Vijay John, Seichi Mita\IEEEauthorrefmark{1}}
	\IEEEauthorblockA{Smart Vehicles Research Center,
		Toyota Technological University\\
		Nagoya, Aichi-Japan\\
		Email: \IEEEauthorrefmark{1}\{ali-boyali, lyuzhemming, vijayjohn, smita\}@toyota-ti.ac.jp,
	           }
     
	}

\maketitle

\begin{abstract}
 In this study, we detail the procedures for designing gain scheduling   controllers by Linear Quadratic $H_\infty$ robust optimization methods in Linear Matrix Inequalities (LMI) framework. The controllers are aimed for steering control of autonomous vehicle. We first construct the Linear Parameter Varying (LPV) vehicle models and synthesize the robust controllers with uncertainty and nominal plants. We choose static output and state feedback controller structure to avoid higher order controllers considering implementation issues. The robust control problems are solved by using different LMI formulations and optimization weights with and without eigenvalue location constraints. The results are compared.
  
\end{abstract}
 
\IEEEpeerreviewmaketitle

\section{Introduction}
 
Self-driving robotic vehicles have been set to become standard in the very near future. We have been witnessing an ongoing race to launch the commercial self-driving vehicles on roads. The notion of self-driving autonomous vehicles not only led to many ethical and legislative challenges but also provided new opportunities to handle the current technical knowledge with the newly emerged requirements. 

The vehicle steering control has been broadly studied. In the literature there are wide variety of vehicle dynamics lateral controllers proposed for steering which meet many performance criteria. However, the maximum performance we can get from the controllers are limited with the parameter uncertainty in the vehicle models, varying operating conditions and available information to the controller.
 
In this paper, we propose robust self-scheduling preview controller syntheses for path tracking. The preview control differs from the Model Predictive Control (MPC) in many aspects. In the preliminary reviews, the reviewers confused the predictive and preview controllers. In the MPC methods, the system model is used to propagate the state space equations through out a prediction horizon predicting the future evolution of the states, and the resulting matrix system is used to solve optimizing input sequence for a given control horizon. The resulting optimization problem is either solved real time. In the explicit MPC problems, the optimization problem solved for operating point trajectories and the solution is stored in look-up tables which are then used during operation instead of real time optimization. However, in the preview control methods, there is no need neither to solve a real time optimization problem nor storing the solutions in look up tables.  It is shown that, the preview control and MPC gives identical results \cite{cole2006predictive} for modeling driver steering control applications. Therefore, the preview controller are more easily realizable then the MPC controller for real time steering control applications. 

The history of preview control dates back to 60's and a plethora of literature on the topic is available. The preview control method we adapted is based on the optimal regulation and servo tracking problem with known input sequence which are theoretically discussed in the optimal control books \cite{anderson1971linear,bitmead1990adaptive}. First applications of the preview control to the path tracking automated driving problems were reported in  \cite{lee1992preview, sharp2006optimal} for 4WD and front wheel steering passenger cars and the time for effective preview, comparison of the controller with the human steering behaviour were discussed. These studies were designed for only constant vehicle speeds and single operating conditions. The uncertainty on the models as well as the parameter variations were not taken into account. 

This paper extends the preview steering controllers by including both parameter variation and uncertainty in the controller design. We contribute to the literature in preview steering control and automated driving with the following solutions;

\begin{itemize}
	\item All the previous steering controllers were proposed for single vehicle speed operating point. We introduce Linear Parameter Varying (LPV) modeling approach and gain scheduling preview steering control by this paper. 
	\item There is no robust control synthesis for preview steering path tracking problems. In this paper we formulate the problem as robust control problem and solve in Linear Matrix Inequalities (LMI) framework.
	\item There is only one LPV-LMI preview control application applied to servomechanism design \cite{takaba2000robust} in which the fast trasient dynamics and pole constraints of the controller are not formulated. In this study, we also extend their approach by adding pole contraints to the LMI formulation which results smooth transient modes. 
	\item  and finally, we give a framework that is flexible to extend the LPV - LMI formulation to beyond the LQ $H_\infty$ using other LMI conditions with dynamic weighing filters.
\end{itemize}

The rest of the paper is organized as following. In Section II, we give brief background information about preview control and its application to vehicle steering and path tracking problem. The LQ $H_{\infty}$ synthesis and robust control of LPV vehicle system are detailed in Section III. In Section IV, the simulation results and comparison are discussed. The paper ends with the conclusion.   
     
\section{Optimal Preview Control}
\subsection{Linear Quadratic (LQ) Preview Control}
 
\begin{figure}[h!]
	\centering
	\includegraphics[width=0.8\linewidth]{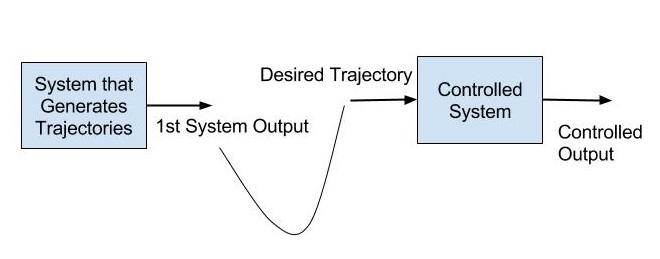}
	\caption{Preview Control Problem}
	\label{fig:prwProblem}
\end{figure}

The reference tracking problem in control theory deals with how to make a system to follow a desired trajectory. If the reference is an output of another system which is known in advance, the problem can be cast into preview control problem (Figure \ref{fig:prwProblem}) and solved by optimization methods. 

The solution to preview control is formulated as a LQ optimization problem in the books \cite{anderson1971linear,bitmead1990adaptive}. The finite and infinite time (time invariant) solutions are obtained by choosing an appropriate cost function and weights (Q for states and R control effort) relating the two systems by state augmentation.

The preview LQ solution was first employed for a vehicle suspension system control and described in  \cite{tomizuka1976optimum, louam1988optimal, prokop1995performance} then the applications extended to the motorcycle and car steering control  problems \cite{sharp2001stability, sharp2001optimal, sharp2006optimal}. In this paper we closely follow \cite{sharp2001optimal} for vehicle-road model development.     

In \cite{sharp2001optimal}, the road is taken as the path generating plant (\ref{Eq:dsystemsr}). The vehicle (\ref{Eq:dsystemsv}) follow the output of this plant. 
\begin{equation}  \scalebox{0.8}{$ 
\Sigma_g=
\begin{cases}
x_r(k+1)=A_{r}x_r(k)+B_{r}y_{ri}(k)\\
y_r(k)=C_{r}x_r(k)
\end{cases} $}
\label{Eq:dsystemsr}
\end{equation}

\begin{equation}  \scalebox{0.8}{$ 
\Sigma_t=
\begin{cases}
x_v(k+1)=A_{v}x_v(k)+B_{v}u_{v}(k)\\
y_v(k)=C_{v}x_v(k)+D_{v}u_v(k)
\end{cases} $}
\label{Eq:dsystemsv}
\end{equation} 
 
In Equation (\ref{Eq:dsystemsr}), $A_r$, $B_r$ and $C_r$ are given by Equation (\ref{Eq:systemMAtrices}) and $y_{ri}$ is the single reference input enters to the system at the preview horizon. The variable $y_{ri}$ is the road state and it represents the deviation of the path preview point at the horizon in the vehicle coordinate system.  

\begin{equation} \scalebox{0.8}{$ 
	\begin{array}{lll}
	A_g=\begin{bmatrix}
	0&1&0&\ldots&0\\
	0&0&1&\ldots&0\\
	\vdots&\vdots&&\ddots&\\
	0&0&0&0&1\\
	0&0&0&0&0\\
	\end{bmatrix}_{NxN}, &   	
	B_g=\begin{bmatrix}
	0\\
	0\\
	\vdots\\
	0\\
	1\\
	\end{bmatrix}_{Nx1}\\
	\\
	C_g=\begin{bmatrix}
	1&0&0&\ldots&0
	\end{bmatrix} &  
	\end{array}$}
\label{Eq:systemMAtrices}
\end{equation}

\begin{figure}[h!]
	\centering
	\includegraphics[width=0.8\linewidth]{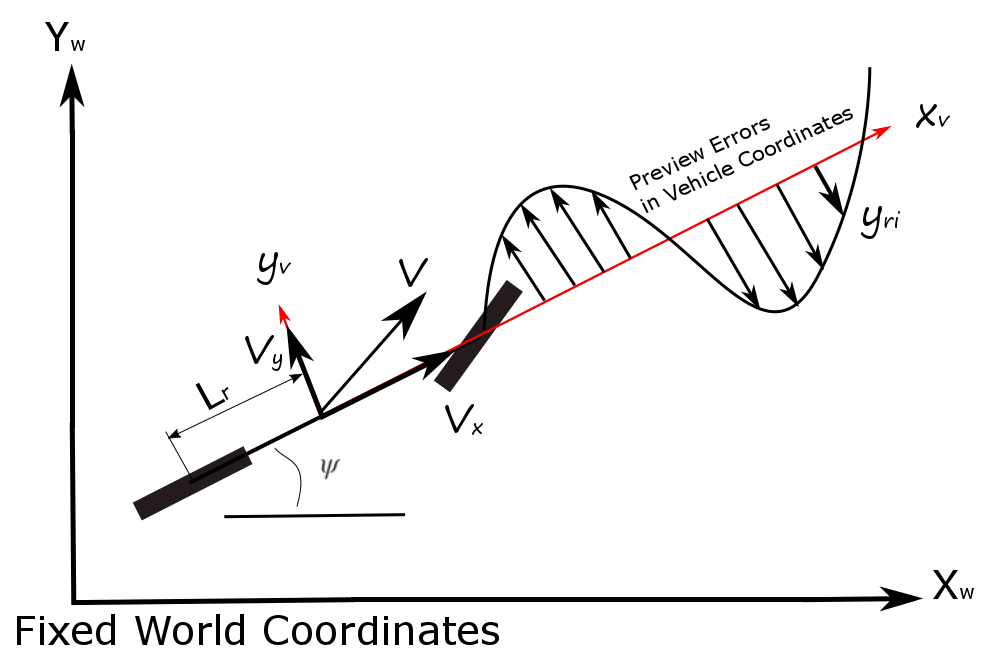}
	\caption{Preview Control Problem}
	\label{fig:singleTrack}
\end{figure}

The vehicle model equations (\ref{Eq:dsystemsv}) are derived from the single track vehicle model for the lateral motion (Figure \ref{fig:singleTrack}). The vehicle local ($x_v$, $y_v$) and the global coordinate frames ($X_w$, $Y_w$) are shown in the figure. The vehicle local longitudinal and lateral speed are represented by the variables $V_x$ and $V_y$ respectively.  

The path tracking error is computed with respect to the vehicle local or global coordinates systems, therefore the appropriate single track model defined in the global and local coordinates should be used for the augmented system. We will give only the error model in local coordinate system for brevity in this paper. 

The vehicle error model is given as described in \cite{rajamani2011vehicle};
\begin{equation} \scalebox{0.8}{$ 
\Sigma_t=
\begin{cases}
x_{ve}(k+1)=A_{ve}x_{ve}(k)+B_{v}u_{v}(k)+B_{\dot{\psi}}{\dot{\psi}}_{des}(k)\\
y_{ve}(k)=C_{ve}x_{ve}(k)+D_{ve}u_v(k)
\end{cases} $}
\label{Eq:derrorv}
\end{equation} 
where the error state transition $A_{ve}$, disturbance and steering input matrices $B_{v}$ and $B_{\dot{\psi}}$ are;
\begin{equation}\scalebox{0.8}{$ 
	A_{ve}=\begin{bmatrix}
	0&1&0&0\\
	0&-\frac{2(C_{{\alpha}f}+C_{{\alpha}r})}{mV_x} &a_{23}&-V_x-\frac{2(C_{{\alpha}f}L_f-C_{{\alpha}r}L_r)}{mV_x} \\
	0&0&0&1\\
	0&-\frac{2(C_{{\alpha}f}L_f-C_{{\alpha}r}L_r)}{I_zV_x} &a_{43}&-\frac{2(C_{{\alpha}f}L_f^2+C_{{\alpha}r}L_r^2)}{I_zV_x} 
	\end{bmatrix}  
	$}
\label{Eq:Ave}
\end{equation} 

\begin{equation}\scalebox{0.8}{$
B_{\dot{\psi}}=\begin{bmatrix}
0\\
-V_x-\frac{2(C_{{\alpha}f}L_f-C_{{\alpha}r}L_r)}{mV_x}\\
0\\
-\frac{2(C_{{\alpha}f}L_f^2+C_{{\alpha}r}L_r^2)}{I_zV_x} 
\end{bmatrix}   
\\and \quad
B_v=\begin{bmatrix}
0\\
\frac{2C_{{\alpha}f}}{m}\\
0\\
\frac{2C_{{\alpha}f}L_f}{I_z}
\end{bmatrix}  $}
\label{Eq:Bv}
\end{equation}

Here,   
\begin{align*}
	 D_{ve}=[0]_{4x2}&\:&
	 \scalebox{1}{$a_{23}=\frac{2(C_{{\alpha}f}+C_{{\alpha}r})}{m}$} &\:& \scalebox{1}{$a_{23}=\frac{2(C_{{\alpha}f}L_f-C_{{\alpha}r}L_r)}{I_z}$}	
\end{align*}

In the equations, $C_{{\alpha}(f,r)}$ are the front and rear cornering stiffness of the tires, $L_{f,r}$ are the distance of the center of gravity of the car to the front and rear axles and $[m, I_z]$ are the mass and moment of inertia of the vehicle around vehicle's $z$ axis respectively. 
The error states from the desired motion $\scalebox{0.9}{$x_v=\begin{bmatrix} y,& V_y, &\Psi, &\dot{\Psi} \end{bmatrix}^T$}$ are the lateral displacement and velocity, heading (yaw) angle and rate in vehicle coordinate system. The road states ($x_r$) contains only the lateral deviation values in the vehicle's current local coordinate frame.  The two systems can now be state augmented in (\ref{Eq:augmentedSys}) to develop an optimal control law by the classical LQR solutions.
\begin{equation}
\scalebox{0.8}{$\begin{bmatrix}
x_{ve}(k+1)\\
x_r(k+1) 
\end{bmatrix}=	\begin{bmatrix}
A_{ve}& 0\\
0 & A_r
\end{bmatrix} \begin{bmatrix}
x_{ve}(k)\\
x_r(k) 
\end{bmatrix}+\begin{bmatrix}
B_{v}\\
0
\end{bmatrix}u_v(k)+\begin{bmatrix}
0\\
B_{r}
\end{bmatrix}y_{ri}
$}
\label{Eq:augmentedSys}																	 
\end{equation}
In the state space equations (\ref{Eq:augmentedSys}), there is no connection between the road and vehicle states. The augmented discrete states are connected with the new measurement matrix $C_{aug}$ to obtain the cost function to be optimized.   
\begin{equation}\scalebox{0.8}{$
\begin{bmatrix*}
x_{aug}(k+1)\\
y_{aug}(k+1) 
\end{bmatrix*}=\begin{bmatrix*}[l]
A_{aug}x_{aug}(k)+B_{r.aug}y_{ri}+B_{v.aug}u_{v}(k)\\
C_{aug}x_{aug}(k)+ D_{v.aug}u_{v}(k)
\end{bmatrix*}$}   
\label{Eq:augmentedSys_Z}																 
\end{equation}
where $\scalebox{0.9}{$x_{aug}(k)=\begin{bmatrix} x_{ve}(k)& y_r(k) \end{bmatrix}^T$}$. 

The road plant provides next preview point and current vehicle heading information to the vehicle at the current time step. Thus the deviation from the next target point  ${\Delta}y=y(k)-y_{r1}$ and current heading angle ${\Psi}{\cong}{\frac{{\Delta}y}{V_xT}} $ are used for the measurement matrix in the augmented error model. In the equation $T$ is the sampling period. In this case the measurement matrix $C_{aug}$ is expressed as;
$$\scalebox{0.8}{$C_{aug}=\begin{bmatrix} 1& 0 &0& 0& -1& 0 & 0&\ldots&0\\0& 0& 1& 0&\frac{1}{V_xT}& -\frac{1}{V_xT}&0&\ldots&0 \end{bmatrix}$}$$

The discrete time cost function (Equation \ref{Eq:costFunc2}) is constructed based on the measurement matrix with the appropriate weights on the outputs. 
\begin{equation}\scalebox{0.8}{$ 
J(x(t_{0}), u(.), t_f))= \sum_{k=0}^{k=N} [(x_{aug})^{T}C_{aug}^{T}QC_{aug}(x_{aug})+u_v^TRu_v]dt $}  
\label{Eq:costFunc2}
\end{equation}
Here (\ref{Eq:costFunc2}) $Q$ is a positive semi-definite matrix containing weights on the diagonal for each measured states while the scaler $R$ is a positive scalar for the control effort. If the assumptions that $[A_{aug}, B_{v}]$ are stabilisable and  $[A_{aug}, C_{aug}]$  detectable hold, the optimal input is computed as $u^*(k)=-Kx_{aug}$. In this case, the state feedback and feedforward controller coefficients are obtained from; 
\begin{equation}\scalebox{0.9}{$
	K=[K_v, K_r]=(R+B_{v.aug}^{T}PB_{v.aug})^{-1}B_{v.aug}^TPA_{aug}$}
	\label{Eq:KvKr}
\end{equation}

The gains $[K_v, K_r]$ are the state feedback and feedforward control coefficients for the vehicle and road states respectively and the matrix $P$ satisfies the Discrete Algebraic Riccati Equation (DARE). The analytic computation of the finite time solutions can be found in the references  \cite{louam1988optimal, prokop1995performance} in detail.  

The preview horizon at which the feed-forward coefficients go zero, depends on the optimization state $Q$ and control weighting $R$ matrices. These weights are adjusted according to the desired response. Although the preview LQ control methods provide many advantage in tracking control, they cannot be applied verbatim for LPV systems. In the following sections, after a brief review of the robust control of LPV plants, we elaborate the LQ $H_{\infty}$ synthesis with gain scheduling.

\subsection{Robust Preview Controllers, Previous Studies}

The advantage of the preview controllers sparkled many study for application of the method in generalized settings for mixed stochastic $LQ{ \mathbin{/}}H_2$ and $LQ{\mathbin{/}}H_\infty$ controllers \cite{hazell2008discrete, takaba2003tutorial, cunha2003terrain, saleh2012optimal}. Almost all of the proposed controllers are based on the partitioning Riccati solution then obtaining the $H_2$ and $H_\infty$ controllers from the Hamiltonian matrix analytically \cite{mianzo2007output, moran2014design, saleh2012optimal}. A few studies propose classical gain scheduling with preview control in which the controller coefficients are computed at the frozen operating grid points. In \cite{thommyppillai2009car}, the authors use the tire parameters as the scheduling parameter. 

A frequency shaped LQ Preview controller is given in \cite{peng1993preview} with the same scheduling variable. The author in \cite{takaba2000robust} formulated polytopic LQ preview control with integral action not for vehicle motion but a servomechanism. The controller coefficients for both preview and the tracking plants are computed in LMI framework. In this study, there is no pole constraint to prevent fast transient dynamics. There are other studies that use the LMI framework in preview control applications, however in these studies, only the feed-forward part of the controllers are computed by LMI methods \cite{silvestrepath, paulino2006affine}. 

In this paper, we formulated the robust steering controller in a similar manner given in \cite{takaba2000robust} for servo-mechanism controller design. In addition to formulating the controllers and solving the problem in LMI framework, we included pole placement constraints to assign the poles of the closed loop system matrix in the desired LMI region. The longitudinal speed of the vehicle is taken as the scheduling variable.

\section{Speed Scheduled Robust Controllers}
In the last decades, the LMI based controller design became a community standard as it allows formulating diverse control problems in a flexible manner. In this paper we formulated the preview control problem in classical robust control theory and solved the resulting matrix inequalities with a highly efficient solver  (MOESP \cite{mosek}) in YALMIP framework \cite{lofberg2004yalmip}  

\begin{figure}[h!]
	\centering
	\includegraphics[scale=0.4]{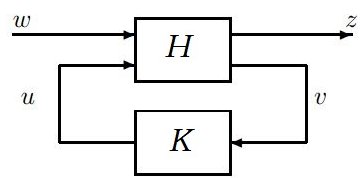}
	\caption{Robust Control Generalized Plant}
	\label{fig:genRobust}
\end{figure}

In robust controller design, the main objective is to reduce the magnitude of the system norms for the specific input-output channels. The general configuration of the system is given in Figure \ref{fig:genRobust}. 
The figure shows a generalized plant $(H)$ with a controller $(K)$ which are connected via a feedback loop. The signals $[\omega, z]$ represent all the exogenous inputs including the uncertainty channels and performance outputs. In preview steering control, the disturbance inputs are the road curvature and the preview input. The uncertainty inputs also enter from this channel to the system when the uncertain state space models are used. The transfer function of this input-output pair $T_{{\omega}z}$ is minimized while satisfying the stability conditions. The state (or static output) feedback controller produce the feedback control signal ($u$) from the measured output ($v$) of the plant. The robust control norm minimization can be computed using LMI convex optimization methods by the Lyapunov stability theory. In this context, a closed loop discrete time system (\ref{Eq:closedLoop2})  is asymptotically stable if and only if there exists a positive symmetric matrix $\mathcal{P}$ satisfying the LMI given by Equation (\ref{Eq:lmi1}). 

\noindent\begin{minipage}{.5\linewidth}
	\begin{align}
		\scalebox{0.8}{$x(k+1)=\mathcal{A}x+\mathcal{B}w$}    \nonumber\\ 
		\scalebox{0.8}{$z(k)=\mathcal{C}x+\mathcal{D}w $} 
		\label{Eq:closedLoop2}  
	\end{align}
\end{minipage}%
\begin{minipage}{.5\linewidth}
	\begin{equation} \scalebox{0.8}{$
		\begin{bmatrix}
			\mathcal{{P}}&\mathcal{AP}\\
			\mathcal{PA^T}&\mathcal{{P}}
		\end{bmatrix}>0$}
	\label{Eq:lmi1}
	\end{equation}
\end{minipage}
 
Similarly, the parameter dependent system models and Lyapunov matrix solution are defined for the LPV systems as;

\noindent\begin{minipage}{.5\linewidth}
	\begin{align}
		\scalebox{0.8}{$x(k+1)=\mathcal{A}(\rho)x+\mathcal{B}(\rho)w$}  \nonumber\\ 
		\scalebox{0.8}{$z(k)=\mathcal{C}(\rho)x+\mathcal{D}(\rho)w $}\nonumber \\ 
	\label{Eq:LPVmodelClosed}  
	\end{align}
\end{minipage}%
\begin{minipage}{.5\linewidth}
\begin{equation}\scalebox{0.8}{$
	\begin{bmatrix}
	\mathcal{{P(\rho)}}&\mathcal{A(\rho)P(\rho)}\\
	\mathcal{P(\rho)A(\rho)^T}&\mathcal{{P(\rho)}}
	\end{bmatrix}>0$}
\label{Eq:lmi2}
\end{equation}
\end{minipage}

where the script letters in (\ref{Eq:LPVmodelClosed}) represent the closed loop system matrices (\ref{Eq:sysmatcl}) and $\rho$ is the varying parameters.

\begin{equation}\scalebox{0.8}{$
\begin{bmatrix}
\mathcal{A}:=A(\rho)+B(\rho)K(\rho)  & \mathcal{B}:=B_w (\rho) \\ 
\mathcal{C}:=C_z(\rho)+D_{zu}(\rho)K(\rho)  & \mathcal{D}:=D_{zw}(\rho) \\ 
\end{bmatrix}$}
\label{Eq:sysmatcl}  
\end{equation} 

In polytopic system model approach, parameter variations are bounded by polytopes with a number of vertices. Therefore, the model for any value of the varying parameters is obtained as a convex combination of the linear models computed at the polytope vertices $\omega_i$ \cite{sename2013robust}. The polytopic representation of a system is given as;

\begin{equation}\scalebox{0.8}
{$
	\Sigma_{LPV}=
	\emph{Co}\begin{Bmatrix}
	\sum_{i=1}^{p} \alpha_{i}\begin{bmatrix}
	\mathcal{A}(\omega_i)&\mathcal{B}(\omega_i)\\
	\mathcal{C}(\omega_i)&\mathcal{D}(\omega_i)
	\end{bmatrix}  
	\end{Bmatrix}
	$}
\end{equation}

where $\sum_{i=1}^{p}\alpha_i=1$ and $\alpha_i>0$, $\alpha_i$ are the barycentric coordinates in the polytope and computed by linear programming or least square methods. The controllers are computed from the closed loop solution as shown in the next subsection.  

\subsection{Robust Preview Steering Control Gain Scheduling}

The vehicle longitudinal speed is measured real-time and is chosen as the scheduling parameter. The speed dependent polytopic models are then constructed by choosing a speed interval. In this application, the speed varies between 3-30 m/s. The lower limit is fixed at 3 m/s is due to the preview length. At the lower speeds the preview length becomes shorter as it depends on the vehicle speed. We use 50 Hz sampling frequency for discrete system and one second preview horizon. At lower speeds such for 1-2 m/s the preview distance corresponds to 1-2 meter which does not contribute much to the controller performance at lower speeds while reducing the controller performance at higher speeds. The real time controller coefficients are computed as the linear combination of all vertice point controllers. 

\begin{figure}[h!]
	\centering
	\includegraphics[scale=0.5]{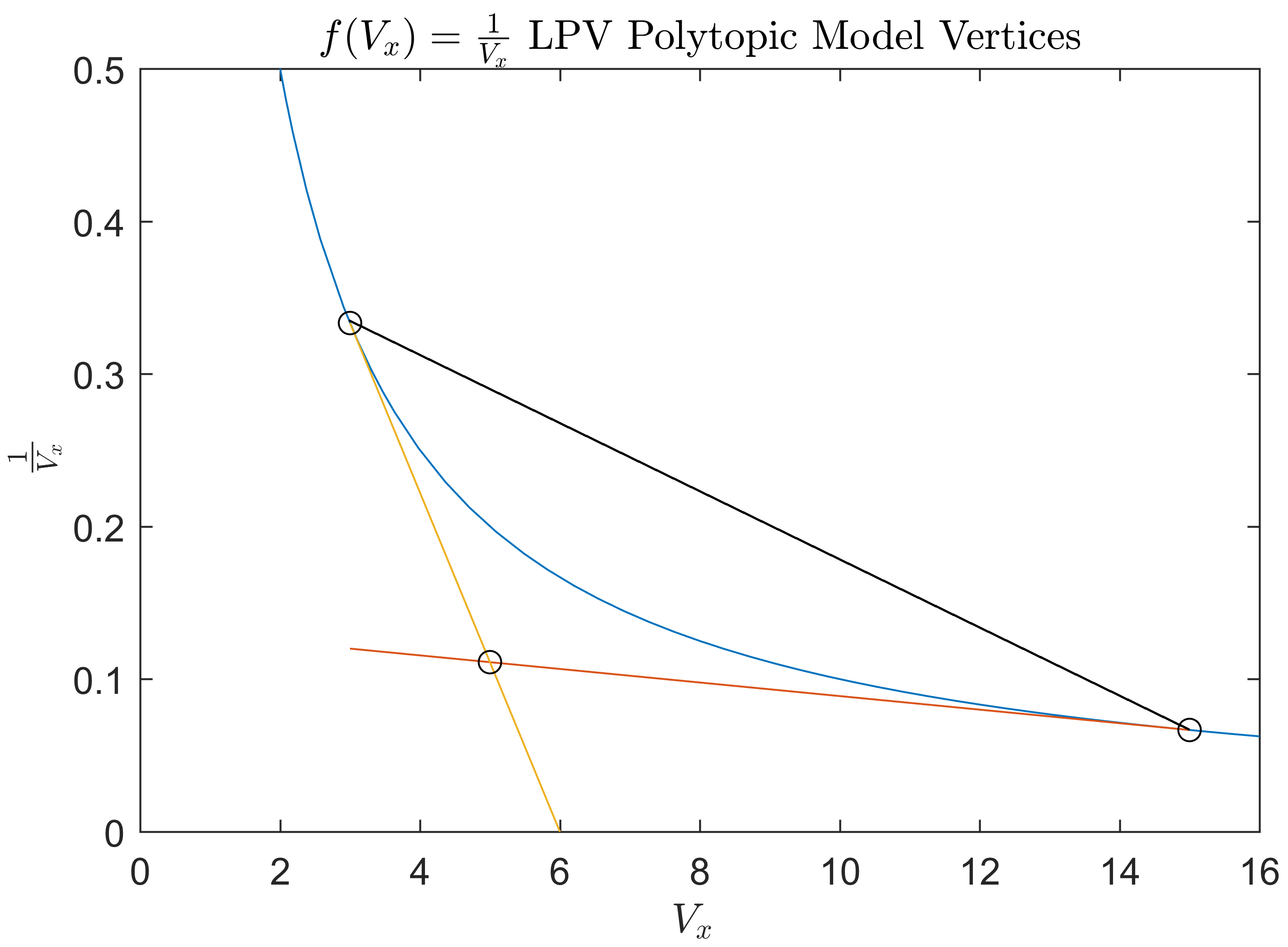}
	\caption{Model Polytope Vertices}
	\label{fig:spdPolytope}
\end{figure}

Upper value of the speed interval can be changed at a cost of higher feedback coefficients. We have two varying variables both are the function of vehicle longitudinal speed ($V_x$ and $1/V_x$). Since both of the varying parameters is a function of vehicle speed, a reduced polytope with three vertices can be constructed as shown in Fig. (\ref{fig:spdPolytope}).

The first two vertex points are the boundary points of the vehicle speed. The third vertex is obtained by taking the intersection point of the lines that are tangent to the speed curve given in the figure. 
 
\begin{figure}[h!]
	\centering
	\includegraphics[scale=0.45]{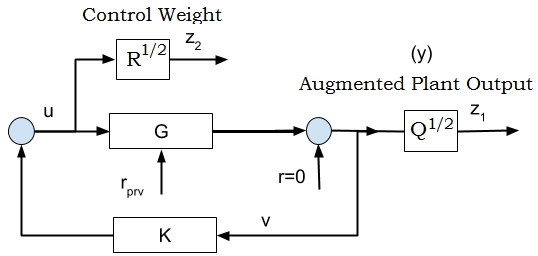}
	\caption{Controller Synthesis - Performance Outputs with Static Weights (LQ)}
	\label{fig:controllerSnth2}
\end{figure}

In robust control, dynamic weighting filters are used to emphasize the desired frequency shapes of the sensitivity ($S$) and complementary sensitivity transfer $T$ functions \cite{zhou1996robust, skogestad2007multivariable}. Control weighting filters are also used to shape the desired control input frequency content. In the LQ $H_\infty$ scheme, the weights on the state and the control signals are constant over all the frequencies. When dynamic weights are used in the optimization, both the number of filters as well as the number of parameters increase. Furthermore, the dynamic filters increase the generalized plant's dimension.The LQ synthesis is more desirable as the weighting parameters are constant at a cost of conservative solution and it is easier to tune. The generalized plant diagram of the LQ settings is given in Fig. (\ref{fig:controllerSnth2}).

The output weighting matrix  $Q=diag(q_1, q_2)$ is a diagonal matrix of two parameters for the lateral deviation of the vehicle from current previewed point and heading angle errors which are defined in the measurement matrix. Two performance channels are represented by $z_{1}$ and $z_{2}$ signals in the diagram. The discrete LMI conditions for $H_\infty$ full state feedback controllers are defined by the following matrix inequities \cite{zhou1995robust,de2002extended}. 

At the polytope vertices ($\omega_i$) where $i=1{\ldots}p$, there exist symmetric positive matrices $P_i$ and  full rank matrix $Z_{i}\in\mathbb{R}^{1,n}$ where $n$ is the total number of states, the following LMI conditions hold for the norm condition of a discrete transfer function from all exogenous inputs to the performance outputs ${\Vert}H_{zw}{\Vert}_2^2<\nu$ of the system;

\scalebox{0.9}{\textbf{\emph{$H_\infty$ Guaranteed Stability Conditions}} for $
	{\Vert}H_{zw}{\Vert}_{\infty}^2<\mu$
}

\begin{equation}
\scalebox{0.9}
{$
	\begin{bmatrix}
	P_i&A_iP_i+B_{u(i)}Z_i&B_{wi}&0\\
	(\star)^T&P_i&0&P_i^TC_{zi}^T+Z_i^TD_{zu(i)}^T\\
	(\star)^T&(\star)^T&I&D_{zw(i)}^T\\
	(\star)^T&(\star)^T&(\star)^T&{\mu}I
	\end{bmatrix}>0
	$}
\label{Eq:hinflmis2}  
\end{equation}

Then the controller coefficients are recovered by the barycentric coordinates using;
\begin{align}
\scalebox{0.8}{$\hat{Z(\alpha)}=\sum_{i}^{p}{\alpha_i}Z_i,$} &\quad& \scalebox{0.8}{$\hat{P(\alpha)}=\sum_{i}^{p}{\alpha_i}P_i,$} \quad&\scalebox{0.8}{$K(\alpha)=\hat{Z(\alpha)}\hat{P(\alpha)}^{-1}$}
\label{Eq:cont1}
\end{align}

The LMI based robust controller solutions tend to result in fast response controller in the transient region. Pole assignment or placement LMI conditions are defined in the LMI equations to shape the transient response of the system. In this study, we require the real part of the all vehicle states are greater than a positive constant.  The pole region constraints for the vehicle states are derived by partitioning the Lyapunov matrices. In this case, the partitioned Lyapunov matrices are expressed as;
 
\begin{equation}\scalebox{0.8}
{$P_i=\begin{bmatrix}
	P_{11,i} & P_{12,i} \\ P_{12,i}^T & P_{22,i}
	\end{bmatrix}{\quad}, P_{11,i}{\in}\mathbb{R}^{nv,nv},  {\quad}P_{22,i}{\in}\mathbb{R}^{nr,nr}
	$}
\label{Eq:partitioned}  
\end{equation}

In Equation (\ref{Eq:partitioned}), $n_{v}$ and $n_{r}$ represent the number of vehicle and road preview states respectively. We can define LMI pole region for the vehicle states by the following inequalities \cite{chilali1999robust, miller2013subspace}.

An LMI region can be represented as a subset;   $$\scalebox{0.8}
{$\mathcal{D}=\{z\in\mathbb{C}:\alpha+z{\beta}+{\bar{z}}{\beta^T}\}$}$$ where the real matrices $\alpha$ and $\beta$ are symmetric such that $\alpha=\alpha^{T}$ and $\beta=\beta^{T}$.

Similarly, $\mathcal{D}$-stable LMI regions can be defined to assign the eigenvalues of closed loop system transition matrices ${A_{ve,}}_{i}$ of the vehicle error model with the symmetric Lyapunov matrices  $P_{11,i}>0$ with the matrix valued equations;
$$\scalebox{0.8}
{$
	\mathcal{M}_{\mathcal{D}}=\alpha{\otimes}{P_{11,i}}+{\beta}{\otimes}({P_{11,i}{A_{ve,}}_{i}})+ {\beta^T}{\otimes}({P_{11,i}{A_{ve,}}_{i}})^T$}$$
   
The positivity constraints for discrete systems can be represented by the region;
$$\scalebox{0.8}
{$\mathcal{D}_p=\{z\in\mathbb{C}:Re(z){\geq}\zeta_{p} \quad \zeta_{p}{\geq}0 \}$}$$   

which corresponds to the following matrix valued function $f_{\mathcal{D}_p}(z)>0$ with

\begin{equation}
\scalebox{0.8}
{$
	 f_{\mathcal{D}_p}(z)=\zeta_{p}\begin{bmatrix}
	 						2&0\\0&-2
	                      \end{bmatrix} + \begin{bmatrix}
	                      0&0\\0&1
	                      \end{bmatrix}z+\begin{bmatrix}
	                      0&0\\0&1
	                      \end{bmatrix}{\hat{z}}
$}	 
\end{equation} 

By these constraints, the eigenvalues of discrete closed loop transition matrix are assigned to have a value greater than a positive constant while the Schur stability conditions (eigenvalues of the closed loop transition matrix is less than one) are guaranteed by the LMI equations (\ref{Eq:hinflmis2}). The closeness of the eigenvalues to the origin of the unit disk defines the speed of transient response. The eigenvalues close the origin yield faster response. 

\section{Results and Discussions}

We computed preview steering state feedback controllers for a mid-range passenger car. The single track nominal parameters are known. We computed two sets of controllers by constructing certain and uncertain vehicle models. In both of the representations, the longitudinal speed is the time varying parameter. The polytopic LPV models are obtained by using longitudinal speed. The state feedback controllers are obtained for different speed intervals such as for the speed intervals [3-15 m/s] and [3-30 m/s]. In the uncertain vehicle models, the tire cornering coefficients are taken as the uncertain parameters. The state space uncertainty is used in the models \cite{zhou1996robust}.  This is easy to do in Matlab by lftdata command. The state feedback coefficients are computed with the certain and uncertain representations in LQ $H_\infty$ solution. However for comparison, we give some results for the models where uncertainty is discarded and all the parameters are treated to be known by solving the robust optimal control problem with extended discrete LMIs given in \cite{de2002extended}. 

\begin{figure}[h!]
	\centering
	\includegraphics[scale=0.4]{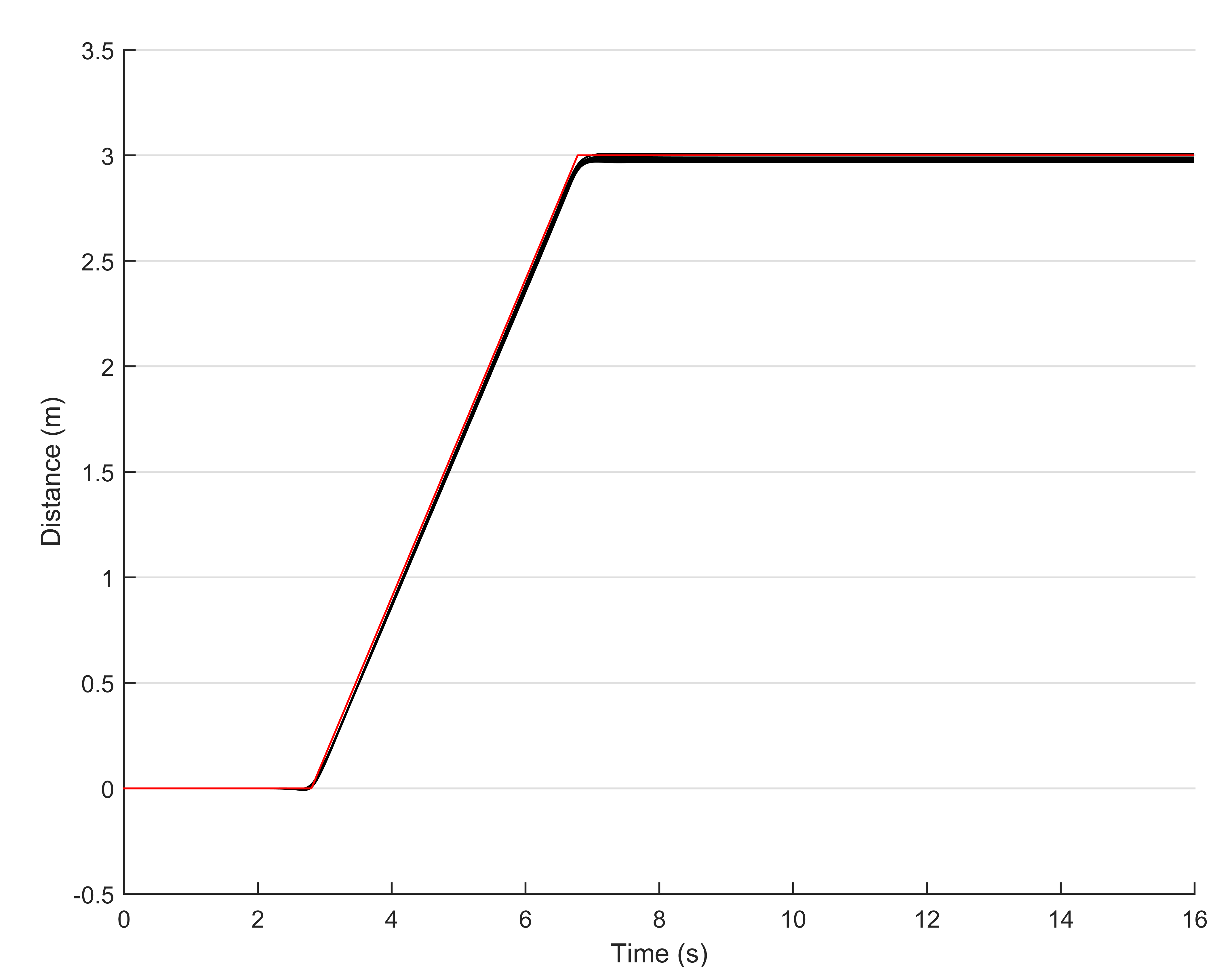}
	\caption{$H_\infty$ Tight Tracking Performance with Dynamic Weights}
	\label{fig:hinftrack}
\end{figure}

We formulated the robust control optimization problem with dynamic weighting filters in mixed sensitivity optimization form with extended LMIs. Three first order weighting filters were defined for sensitivity, complementary sensitivity transfer functions and control signal. We obtained tight tracking performance using $H_\infty$ and $H_2$ norms separately with dynamic filters. However, the computed state feedback coefficients yield fast transient response which is the main characteristics of the LMI solutions as stated in \cite{scherer1997multiobjective}. The LQ $H_\infty$ norm solution also yields the same fast response. The steady state tracking error for all solutions stays in the $\pm$ 3 cm for all speed range. A typical controller with fast response is given for a lane change maneuver in Fig. \ref{fig:hinftrack}. In the figure, the tracking performance for various vehicle speed less than 15 m/s are shown in bold whereas the tracked lane change path is shown by light color line. The vertical and horizontal axes represent the lateral distance and time respectively in the tracking performance figures. 

\begin{figure}[h!]
	\centering
	\includegraphics[scale=0.4]{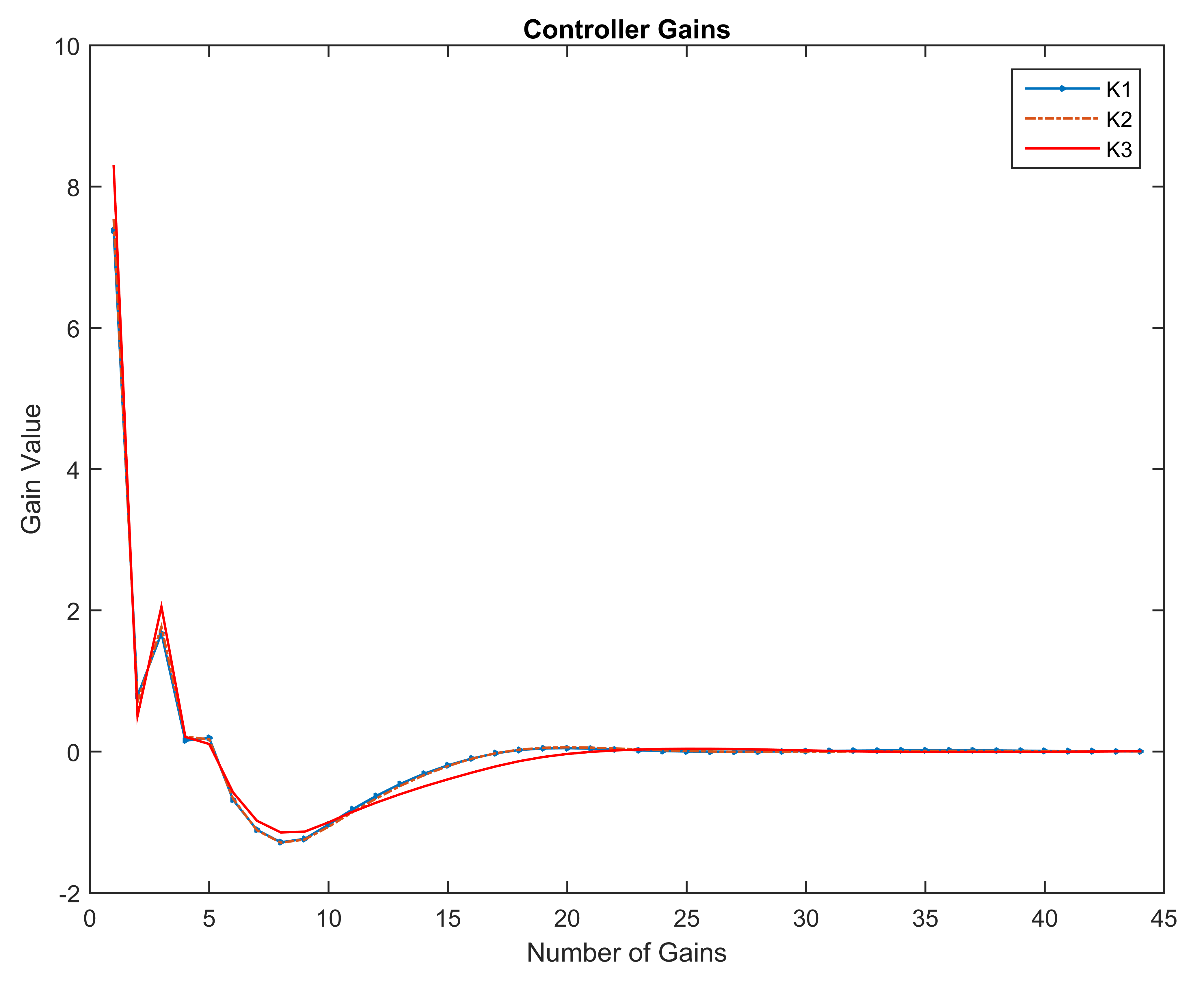}
	\caption{$H_\infty$ Controller Gains computed by using Dynamic Weights}
	\label{fig:hinfgains1}
\end{figure}

The computed state feedback and feedforward controller gains for the vertice models are given in Fig. \ref{fig:hinfgains1} for the solutions with dynamic weighting filters. 

As shown in the figure, the preview controller gain vectors (K1, K2, and K3) at the polytope vertices go zero after the \nth{30} preview point. The first four controller gains in the figure correspond to the feedback gains which are high that may give rise to implementation problems. We also obtain fast response within the LQ scheme. In order to slow the transient response, it is essential to use additional pole or eigenvalue constraints in the solution. We set a positivity constraint to push the four eigenvalues of the vehicle plant transition matrix away from origin by partitioning the Lyapunov matrix. In this way, we can obtain slower response with lower constant optimization weights  by setting    $Q=diag(0.95, \:3e{-3})$ and the control signal weight $R=0.25$. The tracking performance and control gain results of the LQ robust control solution for nominal plant without uncertainty are given in the following figures (\ref{fig:hinfLQtracking}) for the speed interval [3-30 m/s]. On the cornering points, the response gets smoother in Figure \ref{fig:hinfLQtracking}.

\begin{figure}[h!]
	\centering
	\includegraphics[scale=0.4]{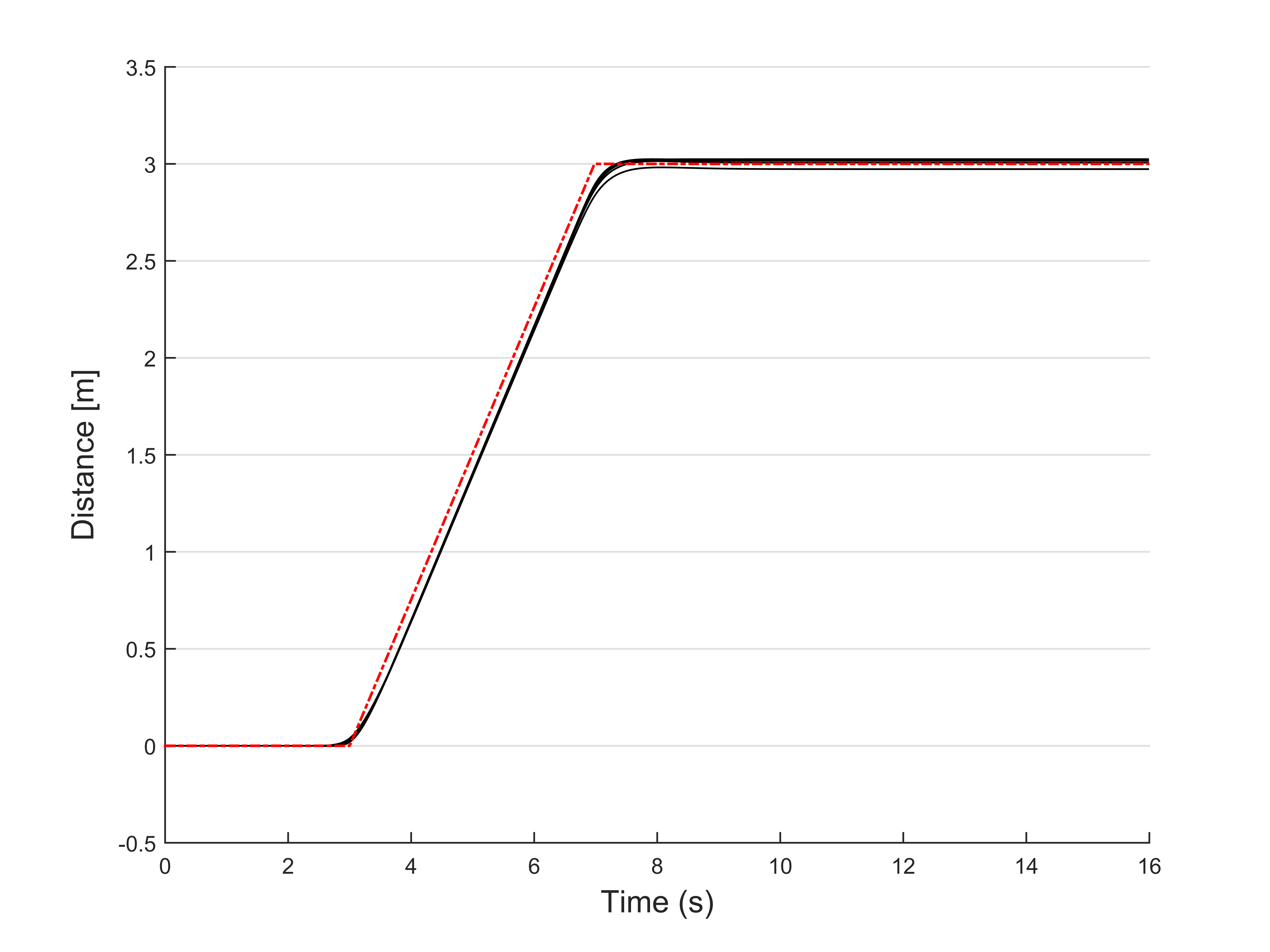}
	\caption{$H_\infty$ LQ Solution Tracking Performance with Constant Weights}
	\label{fig:hinfLQtracking}
\end{figure}  

\begin{figure}[h!]
	\centering
	\includegraphics[scale=0.4]{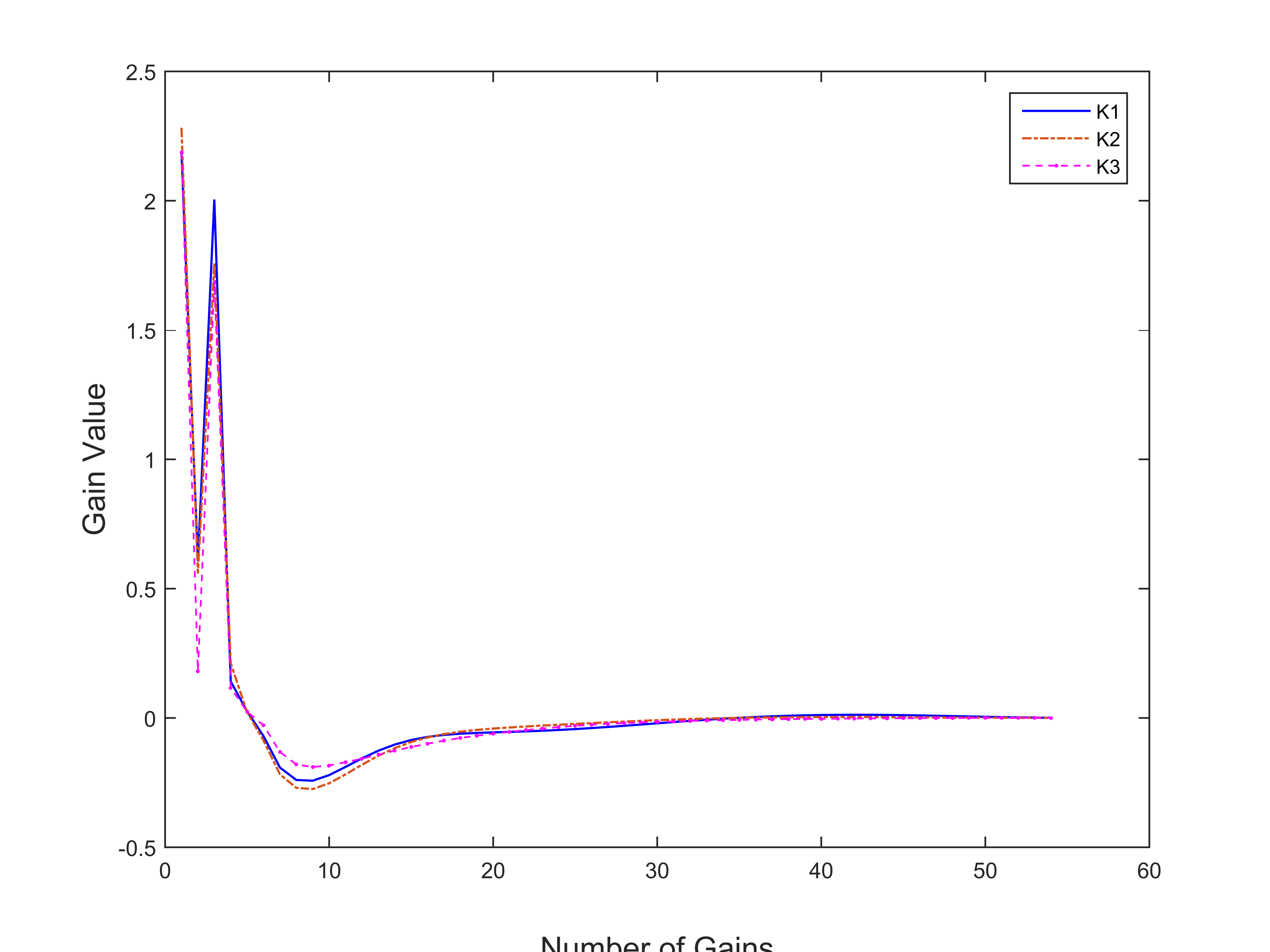}
	\caption{$H_\infty$ LQ solution Controller Gains}
	\label{fig:hinfLQgains}
\end{figure}

The corresponding gain vectors at each of the polytope vertices are shown in Fig. \ref{fig:hinfLQgains}.  We show the gains computed in LQ scheme with uncertainty. As seen in the figure, the level of the both feedback and feedforward gains reduce significantly. In Figure \ref{fig:hinfLQgains}, only preview gains for 50 preview points are shown closely.

\begin{figure}[h!]
	\centering
	\includegraphics[scale=0.4]{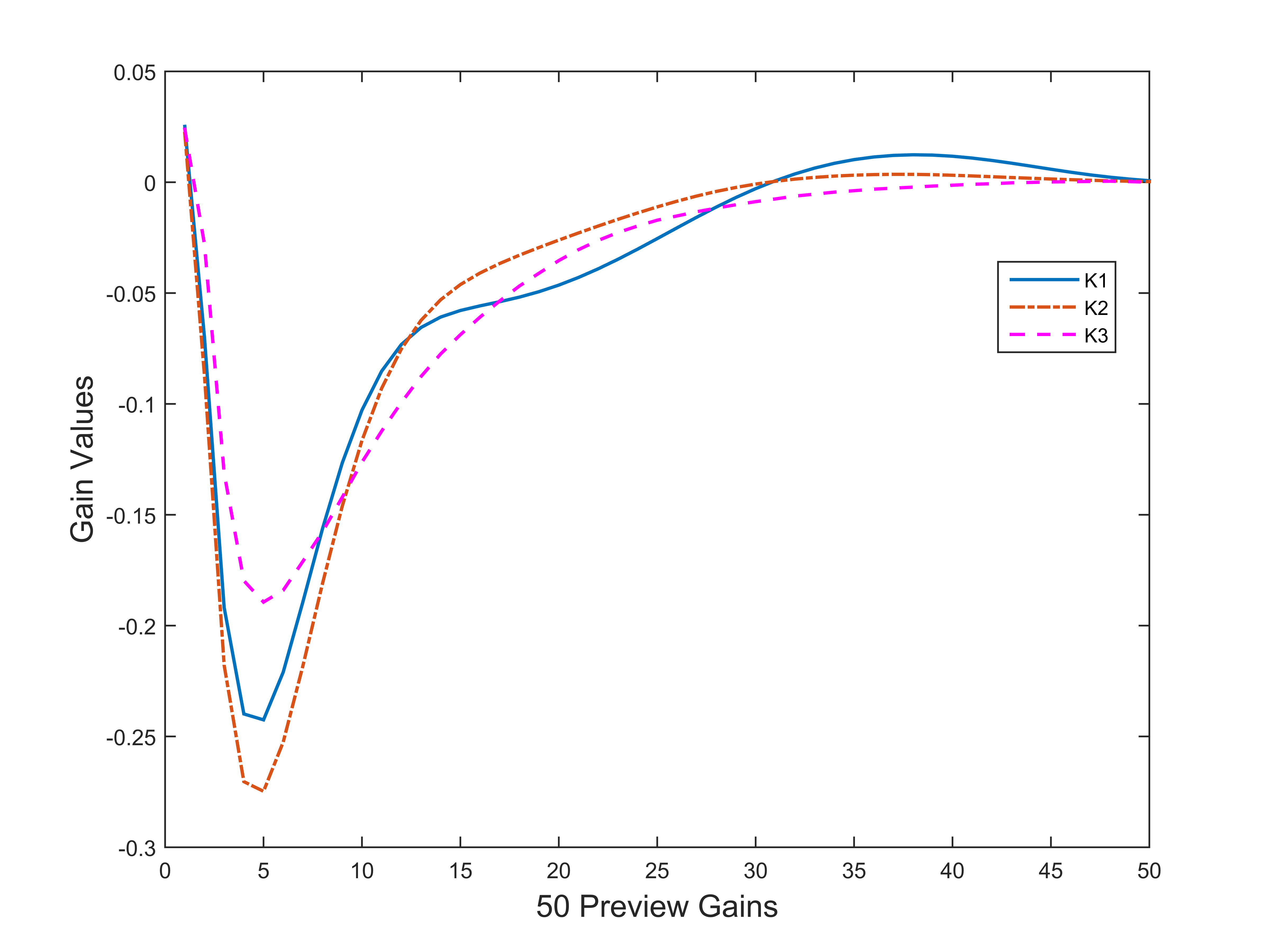}
	\caption{$H_\infty$ LQ solution Controller Preview Gains}
	\label{fig:hinfLQprvgains}
\end{figure}  

\section{Conclusion}
In this study we presented robust preview controllers in LQ $H_\infty$ form for certain and uncertain plants. The robust optimization problems were solved using constant weighting matrices for certain and uncertain model parameters. The tire cornering stiffnesses are treated as the uncertain with a range of $\pm$30 \%. The resulting control and tracking performance were compared to the solutions obtained when dynamic weights are used with extended LMI formulations. The robust optimization solved in LMI framework yielded fast transient response with higher controller gains. We added eigenvalue constraints in the solution to prevent the higher gains and fast response and obtained more implementable controller with lower gains for 50 preview points. The results show that in the current settings 50 preview points are sufficient as the number beyond 50, the preview gains goes zero. The computational time is short and fast, however when the system dimensions increase with the dynamic weighting filters, and with number of optimization variables in the extended LMI formulation the computation times get longer. Considering the tuning process and re-running the scripts make the computations much longer with addition of the pole constraints. We solved the overall problem within LQ scheme and obtained satisfactory solutions in shorter times.

\section*{Acknowledgment}



\bibliographystyle{IEEEtran}
\bibliography{ascc2017}

\end{document}